\documentclass[reqno]{amsart}%
\usepackage{amssymb,mathtools,enumerate}

\usepackage{ifpdf}
\ifpdf
 \usepackage[hyperindex]{hyperref}
\else
 \expandafter\ifx\csname dvipdfm\endcsname\relax
 \usepackage[hypertex,hyperindex]{hyperref}%
 \else
 \usepackage[dvipdfm,hyperindex]{hyperref}%
 \fi
\fi

\allowdisplaybreaks[4]
\numberwithin{equation}{section}

\theoremstyle{plain}
\newtheorem{thm}{Theorem}[section]

\theoremstyle{remark}
\newtheorem{rem}{Remark}[section]

\DeclareMathOperator{\bell}{B}
\DeclareMathOperator{\te}{e}
\DeclareMathOperator{\ti}{i}
\DeclareMathOperator{\td}{d\!}

\begin{document}

\title[Concise proofs of formulas for complete Bell polynomials]
{Concise and elegant proofs of three formulas for complete Bell polynomials}

\author[F. Qi]{Feng Qi*}
\address{17709 Sabal Court, University Village, Dallas, TX 75252-8024, USA}
\email{\href{mailto: F. Qi <qifeng618@gmail.com>}{qifeng618@gmail.com}}
\urladdr{\url{https://orcid.org/0000-0001-6239-2968}}

\begin{abstract}
In the paper, in light of the generating function of the complete Bell polynomials and other techniques, the author presents concise and elegant proofs of three formulas for the complete Bell polynomials.
\end{abstract}



\thanks{*Corresponding author}

\thanks{This paper was typeset using \AmS-\LaTeX}

\maketitle

\section{Introduction}
\subsection{Hoffman's formula}
In 1992, Hoffman~\cite{Hoffman-Pacific-1992} and Zagier~\cite{Don.1992-94} introduced and investigated the multiple harmonic series (currently known as the multiple zeta functions or the multiple zeta values)
\begin{equation*}
\zeta^*(z_1,z_2,\dotsc,z_k)=\sum_{n_1\ge n_2\ge \dotsm\ge n_k\ge1}\frac{1}{n_1^{z_1}n_2^{z_2}\dotsm n_k^{z_k}}
\end{equation*}
and
\begin{equation*}
\zeta(z_1,z_2,\dotsc,z_k)=\sum_{n_1>n_2>\dotsm>n_k\ge1}\frac{1}{n_1^{z_1}n_2^{z_2}\dotsm n_k^{z_k}}
\end{equation*}
for $k\in\mathbb{N}=\{1,2,\dotsc\}$ and $\Re(z_k)>1$. This is a natural generalization of the classical Riemann zeta function 
\begin{equation*}
\zeta(z)=\sum_{n=1}^{\infty}\frac{1}{n^z}, \quad\Re(z)>1.
\end{equation*}
For more information about $\zeta(z)$ and its recent properties, please refer to~\cite[Section~3.5]{Temme-96-book} and the review article~\cite{MIA-9509.tex}.
The multiple zeta function $\zeta(z_1,z_2,\dotsc,z_k)$ plays an important role in quantum physics and in knot theory, please refer to the monographs~\cite{Gil-Fresan-2020, Zhao-zeta-multiple}.
\par
In order to derive
\begin{equation*}
\zeta(\underbrace{2,2,\dotsc,2}_k)=\frac{\pi^{2k}}{(2k+1)!}, \quad k\in\mathbb{N}
\end{equation*}
in~\cite[Corollary~2.3]{Hoffman-Pacific-1992}, Hoffman established in~\cite[Proposition~2.4]{Hoffman-Pacific-1992} the following elegant formula.

\begin{thm}[{\cite[Proposition~2.4]{Hoffman-Pacific-1992}}] \label{Hoffman-form-thm}
For $k\in\mathbb{N}$,
\begin{equation}\label{Hoffman-ID-1992}
\sum_{\substack{\sum_{i=1}^ki \ell_i=k\\ \ell_1,\ell_2,\dotsc,\ell_k\in\mathbb{N}_0}} \prod_{i=1}^k\frac{1}{\ell_i!}  \biggl[\frac{B_{2i}}{(2i)(2i)!}\biggr]^{\ell_i}
=\frac{1}{2^{2k}(2k+1)!},
\end{equation}
where $B_{2i}$ denotes the Bernoulli numbers generated~\cite[p.~3, Eq.~(1.1)]{Temme-96-book} by
\begin{equation*}
\frac{z}{\te^z-1}=\sum_{i=0}^{\infty}B_i\frac{z^i}{i!}=1-\frac{1}{2}z+\sum_{i=1}^{\infty}B_{2i}\frac{z^{2i}}{(2i)!}, \quad |z|<2\pi.
\end{equation*}
\end{thm}

\subsection{Gen\v{c}ev's extensions of Hoffman's formula}
In 2024, three decades later, among other things, Gen\v{c}ev~\cite{Marian-Gencev-Results-2024} extended the formula~\eqref{Hoffman-ID-1992} in Theorem~\ref{Hoffman-form-thm} via specific zeta-like series in terms of the Bernoulli numbers $B_{2k}$, the Euler numbers $E_{2k}$, and the Catalan numbers $C_k$ as follows.

\begin{thm}[{\cite[Theorem~2.1]{Marian-Gencev-Results-2024}}]\label{Marian-Gencev-Results-2024Theorem2.1}
For $k\in\mathbb{N}_0=\{0,1,2,\dotsc\}$, we have
\begin{equation}\label{Gencev-ID1-2024}
\sum_{\substack{\sum_{j=1}^kj \ell_j=k\\ \ell_1,\ell_2,\dotsc,\ell_k\in\mathbb{N}_0}} \frac{1}{\prod_{j=1}^k\ell_j!} \prod_{j=1}^k\biggl[\frac{\epsilon B_{2j}}{(2j)(2j)!}\biggr]^{\ell_j}
=
\begin{dcases}
\frac{1}{2k+1}\frac{1}{(4k)!!}, & \epsilon=1;\\
\frac{2-2^{2k}}{(4k)!!}B_{2k}, & \epsilon=-1.
\end{dcases}
\end{equation}
\end{thm}

\begin{thm}[{\cite[Theorem~2.2]{Marian-Gencev-Results-2024}}]\label{Marian-Gencev-Results-2024Theorem2.2}
For $k\in\mathbb{N}_0$, we have
\begin{equation}\label{Gencev-ID2-2024}
\sum_{\substack{\sum_{j=1}^kj \ell_j=k\\ \ell_1,\ell_2,\dotsc,\ell_k\in\mathbb{N}_0}} \frac{1}{\prod_{j=1}^k\ell_j!} \prod_{j=1}^k\biggl[\bigl(4^j-1\bigr)\frac{\epsilon B_{2j}}{(2j)(2j)!}\biggr]^{\ell_j}
=
\begin{dcases}
\frac{1}{(4k)!!}, & \epsilon=1;\\
\frac{1}{(4k)!!}E_{2k}, & \epsilon=-1,
\end{dcases}
\end{equation}
where $E_{2k}$ denotes the Euler numbers generated~\cite[p.~804, 23.1.2]{abram} by
\begin{equation*}
\frac{2}{\te^z+\te^{-z}}=\sum_{k=0}^{\infty}E_k\frac{z^k}{k!}
=\sum_{k=0}^{\infty}E_{2k}\frac{z^{2k}}{(2k)!}, \quad |z|<\frac\pi2.
\end{equation*}
\end{thm}

\begin{thm}[{\cite[Theorem~2.5]{Marian-Gencev-Results-2024}}]\label{Marian-Gencev-Results-2024Theorem2.5}
For $k\in\mathbb{N}$, we have
\begin{equation}\label{Gencev-ID3-2024}
\sum_{\substack{\sum_{j=1}^kj \ell_j=k\\ \ell_1,\ell_2,\dotsc,\ell_k\in\mathbb{N}_0}} \frac{1}{\prod_{j=1}^k\ell_j!} \prod_{j=1}^k\biggl[\frac{\epsilon}{2j}\binom{2j}{j}\biggr]^{\ell_j}
=\begin{cases}
C_k, & \epsilon=1;\\
-C_{k-1}, & \epsilon=-1,
\end{cases}
\end{equation}
where $C_k=\frac{1}{k+1}\binom{2k}{k}$ for $k\in\mathbb{N}_0$ denotes the Catalan numbers generated~\cite{Koshy-B-2009, AADM-3251.tex, Stanley-Catalan-2015} by
\begin{equation}\label{Catablan-gen-F}
\frac2{1+\sqrt{1-4x}\,}
=\sum_{k=0}^\infty C_k x^k
=1+ x+ 2x^2+ 5x^3+\dotsm, \quad |x|<\frac{1}{4}.
\end{equation}
\end{thm}

\subsection{He--Qi's reformulations of Gen\v{c}ev's formulas}
In~\cite[p.~134, Theorem~A]{Comtet-Combinatorics-74}, the Bell polynomials of the second kind, also known as the partial Bell polynomials, are defined by
\begin{equation}\label{partial-Bell-Polyn-Dfn}
\bell_{k,j}(a_1,a_2,\dotsc,a_{k-j+1})
=\sum_{\substack{\\ \sum_{i=1}^{k-j+1}i\ell_i=k\\
\sum_{i=1}^{k-j+1}\ell_i=j\\ \ell_1,\ell_2,\dotsc,\ell_{k-j+1}\in\mathbb{N}_0}}\frac{k!}{\prod_{i=1}^{k-j+1}\ell_i!} \prod_{i=1}^{k-j+1}\biggl(\frac{a_i}{i!}\biggr)^{\ell_i}
\end{equation}
for $j,k\in\mathbb{N}_0$ satisfying $k\ge j$, with the special cases $\bell_{0,0}(a_1)=1$ and
\begin{equation*}
\bell_{k,0}(a_1,a_2,\dotsc,a_{k+1})=0, \quad k\in\mathbb{N}.
\end{equation*}
In~\cite[p.~412, Definition~11.1]{Charalambides-book-2002} and~\cite[p.~134]{Comtet-Combinatorics-74}, the complete Bell polynomials, denoted by $\bell_k(a_1,a_2,\dotsc,a_k)$, are defined by $\bell_0=1$ and
\begin{equation}\label{complete-Bell-Polyn-Dfn}
\bell_k(a_1,a_2,\dotsc,a_k)=\sum_{j=1}^k\bell_{k,j}(a_1,a_2,\dotsc,a_{k-j+1}), \quad k\in\mathbb{N}.
\end{equation}
For more information on the partial and complete Bell polynomials, please refer to~\cite[Section~3.3, pp.~133--137]{Comtet-Combinatorics-74}.
Directly from the relation~\eqref{complete-Bell-Polyn-Dfn}, it follows that
\begin{equation}\label{A2-Reviewer}
\frac{\bell_k(a_1,a_2,\dotsc,a_k)}{k!}
=\sum_{\substack{\sum_{i=1}^ki\ell_i=k\\ \ell_1,\ell_2,\dotsc,\ell_k\in\mathbb{N}_0}} \prod_{i=1}^k\biggl[\frac{1}{\ell_i!}\biggl(\frac{a_i}{i!}\biggr)^{\ell_i}\biggr], \quad k\in\mathbb{N}.
\end{equation}
The complete Bell polynomials $\bell_k(a_1,a_2,\dotsc,a_k)$ have an exponential generating function
\begin{equation}\label{exp-complete-Bell}
\exp\Biggl(\sum_{k=1}^{\infty}a_k\frac{z^k}{k!}\Biggr)=\sum_{k=0}^{\infty}\bell_k(a_1,a_2,\dotsc,a_k)\frac{z^k}{k!},
\end{equation}
which can be found in~\cite[Chapter~11]{Charalambides-book-2002} and~\cite[p.~134]{Comtet-Combinatorics-74}.
\par
In 2024, observing and comparing the expressions on the left-hand sides of the identities~\eqref{Hoffman-ID-1992} through~\eqref{Gencev-ID3-2024} with the expression on the right-hand side of~\eqref{A2-Reviewer}, He and Qi~\cite{RIMA-D-23-01264.tex} reformulated and alternatively proved the formulas~\eqref{Gencev-ID1-2024}, \eqref{Gencev-ID2-2024}, and~\eqref{Gencev-ID3-2024} in Theorems~\ref{Marian-Gencev-Results-2024Theorem2.1}, \ref{Marian-Gencev-Results-2024Theorem2.2}, and~\ref{Marian-Gencev-Results-2024Theorem2.5} in terms of the complete Bell polynomials $\bell_k$ as follows.

\begin{thm}[{\cite[Theorem~4]{RIMA-D-23-01264.tex}}]\label{Theorem2.1Marian-Gencev-Results-2024-Bell}
For $k\in\mathbb{N}_0$ and $\epsilon=\pm1$, we have
\begin{equation}\label{Gencev-ID1-2024-Bell}
\bell_k\biggl(\frac{\epsilon}{2}\frac{B_{2}}{2!}, \frac{\epsilon}{2}\frac{B_{4}}{4!}, \dotsc,\frac{\epsilon}{2}\frac{(k-1)!B_{2k}}{(2k)!}\biggr)
=
\begin{dcases}
\frac{1}{2k+1}\frac{k!}{(4k)!!}, & \epsilon=1;\\
\frac{k!\bigl(2-2^{2k}\bigr)}{(4k)!!}B_{2k}, & \epsilon=-1.
\end{dcases}
\end{equation}
\end{thm}

\begin{thm}[{\cite[Theorem~5]{RIMA-D-23-01264.tex}}]\label{Theorem2.2Marian-Gencev-Results-2024-Bell}
For $k\in\mathbb{N}_0$ and $\epsilon=\pm1$, we have
\begin{equation}\label{Gencev-ID2-2024-Bell}
\bell_k\biggl(\frac{3\epsilon}{2}\frac{B_{2}}{2!}, \frac{15\epsilon}{2}\frac{B_{4}}{4!}, \dotsc,\frac{\bigl(2^{2k}-1\bigr)\epsilon}{2}\frac{(k-1)!B_{2k}}{(2k)!}\biggr)
=
\begin{dcases}
\frac{k!}{(4k)!!}, & \epsilon=1;\\
\frac{k!}{(4k)!!}E_{2k}, & \epsilon=-1.
\end{dcases}
\end{equation}
\end{thm}

\begin{thm}[{\cite[Theorem~6]{RIMA-D-23-01264.tex}}]\label{Theorem2.5Marian-Gencev-Results-2024-Bell}
For $k\in\mathbb{N}$ and $\epsilon=\pm1$, we have
\begin{equation}\label{Gencev-ID3-2024-Bell}
\bell_k\biggl(\epsilon, 3\epsilon, 20\epsilon, \dotsc,\frac{\epsilon(k-1)!}{2}\binom{2k}{k}\biggr)
=
\begin{dcases}
k!C_k, & \epsilon=1;\\
-k!C_{k-1}, & \epsilon=-1.
\end{dcases}
\end{equation}
\end{thm}

\subsection{He--Qi's generalizations of Gen\v{c}ev's formulas}
In light of the generating function~\eqref{exp-complete-Bell} and other intricate techniques, He and Qi~\cite{RIMA-D-23-01264.tex} not only alternatively proved Theorems~\ref{Theorem2.1Marian-Gencev-Results-2024-Bell} through~\ref{Theorem2.5Marian-Gencev-Results-2024-Bell}, but also generalized Theorems~\ref{Theorem2.1Marian-Gencev-Results-2024-Bell} through~\ref{Theorem2.5Marian-Gencev-Results-2024-Bell} from the specific cases $\epsilon=\pm1$ to the general case $\epsilon\in\mathbb{R}$. We now recite these generalizations as follows.

\begin{thm}[{\cite[Theorem~7]{RIMA-D-23-01264.tex}}]\label{Theorem2.1Marian-Gencev-Results-2024-Bell-gen}
For $k\in\mathbb{N}$ and $\epsilon\in\mathbb{R}$, we have
\begin{multline}\label{Gencev-ID1-2024-Bell-gen}
\bell_k\biggl(\epsilon\frac{B_{2}}{2!}, \epsilon\frac{B_{4}}{4!}, \dotsc,\epsilon\frac{(k-1)!B_{2k}}{(2k)!}\biggr)\\
=\frac{k!}{(2k)!}\sum_{\ell=1}^{2k}\frac{(-2\epsilon)_\ell}{\ell!} \sum_{j=1}^\ell(-1)^j\binom{\ell}{j} \frac{T(2k+j,j)}{\binom{2k+j}{j}},
\end{multline}
where the rising factorial $(z)_\ell$ for $z\in\mathbb{C}$ is defined~\cite[p.~256]{abram} by
\begin{equation}\label{rising-Factorial}
(z)_\ell=\prod_{\ell=0}^{\ell-1}(z+\ell)
=
\begin{cases}
z(z+1)\dotsm(z+\ell-1), & \ell\ge1\\
1, & \ell=0
\end{cases}
\end{equation}
and the central factorial numbers of the second kind $T(p,q)$ for $p,q\in\mathbb{N}_0$ can be computed~\cite[pp.~213--214]{Riordan-B-1968} by
\begin{equation*}
T(p,q)=\frac1{q!} \sum_{k=0}^{q}(-1)^k\binom{q}{k}\biggl(\frac{q}2-k\biggr)^p, \quad p,q\in\mathbb{N}_0
\end{equation*}
with $T(q,q)=1$ for $q\in\mathbb{N}_0$ and $T(p,0)=0$ for $p\in\mathbb{N}$.
\end{thm}

\begin{thm}[{\cite[Theorem~8]{RIMA-D-23-01264.tex}}]\label{Theorem2.2Marian-Gencev-Results-2024-Bell-gen}
For $k\in\mathbb{N}_0$ and $\epsilon\in\mathbb{R}$, we have
\begin{multline}\label{Gencev-ID2-2024-Bell-gen}
\bell_k\biggl(3\epsilon\frac{B_{2}}{2!}, 15\epsilon\frac{B_{4}}{4!}, \dotsc, \bigl(2^{2k}-1\bigr)(k-1)!\epsilon\frac{B_{2k}}{(2k)!}\biggr)\\
=\frac{k!}{2^{2k}(2k)!} \sum_{\ell=0}^{2k} \frac{(-2\epsilon)_\ell}{\ell!} \sum_{m=0}^\ell\frac{(-1)^m}{2^m}\binom{\ell}{m}
\sum_{q=0}^m\binom{m}{q}(2q-m)^{2k},
\end{multline}
where $0^0$ is understood as $1$.
\end{thm}

\begin{thm}[{\cite[Theorem~9]{RIMA-D-23-01264.tex}}]\label{Theorem2.5Marian-Gencev-Results-2024-Bell-gen}
For $k\in\mathbb{N}_0$ and $\epsilon\in\mathbb{R}$, we have
\begin{equation}\label{Gencev-ID3-2024-Bell-gen}
\bell_k\biggl(2\epsilon, 6\epsilon, 40\epsilon, \dotsc,\epsilon(k-1)!\binom{2k}{k}\biggr)
=\sum_{\ell=0}^k(2\epsilon)_{k-\ell} \binom{k+\ell-1}{2\ell} 2^{\ell}(2\ell-1)!!.
\end{equation}
\end{thm}

\subsection{Xu's concise forms for He--Qi's formulas}

In early 2026, using the concept of the complete Bell polynomials $\bell_k(a_1,a_2,\dotsc,a_k)$ and the equation~\eqref{exp-complete-Bell} again, Xu~\cite{Xu-AIMS-Math-2025} presented concise forms for the formulas~\eqref{Gencev-ID1-2024-Bell-gen}, \eqref{Gencev-ID2-2024-Bell-gen}, and~\eqref{Gencev-ID3-2024-Bell-gen} in Theorems~\ref{Theorem2.1Marian-Gencev-Results-2024-Bell-gen}, \ref{Theorem2.2Marian-Gencev-Results-2024-Bell-gen}, and~\ref{Theorem2.5Marian-Gencev-Results-2024-Bell-gen}. We recite his results as follows.

\begin{thm}[{\cite[Theorem~3.1]{Xu-AIMS-Math-2025}}]\label{Xu-thm1}
For $k\in\mathbb{N}_0$ and $\epsilon\in\mathbb{R}$, we have
\begin{equation}\label{Xu-Eq1}
\bell_k\biggl(\epsilon\frac{B_{2}}{2!}, \epsilon\frac{B_{4}}{4!}, \dotsc,\epsilon\frac{(k-1)!B_{2k}}{(2k)!}\biggr)
=\frac{k!}{(2k)!}B_{2k}^{(-2\epsilon)}(-\epsilon),
\end{equation}
where the generalized Bernoulli polynomials $B_{2k}^{(\sigma)}$ for $\sigma\in\mathbb{C}$ are generated~\cite[p.~4]{Temme-96-book} by
\begin{equation}\label{Genaralized-BPolyn}
\biggl(\frac{z}{\te^z-1}\biggr)^\sigma\te^{x z}=\sum_{k=0}^{\infty}B_k^{(\sigma)}(x)\frac{z^k}{k!}, \quad |z|<2\pi.
\end{equation}
\end{thm}

\begin{thm}[{\cite[Corallary~3.1]{Xu-AIMS-Math-2025}}]\label{Xu-thm2}
For $k\in\mathbb{N}_0$ and $\epsilon\in\mathbb{R}$, we have
\begin{equation}\label{Xu-Eq2}
\bell_k\biggl(3\epsilon\frac{B_{2}}{2!}, 15\epsilon\frac{B_{4}}{4!}, \dotsc, \bigl(2^{2k}-1\bigr)\epsilon\frac{(k-1)!B_{2k}}{(2k)!}\biggr) =\frac{k!}{(2k)!}E_{2k}^{(-2\epsilon)}(-\epsilon),
\end{equation}
where the generalized Euler polynomials $E_{2k}^{(\sigma)}$ for $\sigma\in\mathbb{C}$ are generated~\cite[p.~16]{Temme-96-book} by
\begin{equation}\label{Gen-Euler-Polyn}
\biggl(\frac{2}{\te^z+1}\biggr)^\sigma\te^{x z}=\sum_{k=0}^{\infty}E_k^{(\sigma)}(x)\frac{z^k}{k!}, \quad |z|<\pi.
\end{equation}
\end{thm}

\begin{thm}[{\cite[Theorem~3.3]{Xu-AIMS-Math-2025}}]\label{Xu-thm3}
For $k\in\mathbb{N}$ and $\epsilon\in\mathbb{R}$, we have
\begin{equation}\label{Xu-Eq3}
\bell_k\biggl(2\epsilon, 6\epsilon, 40\epsilon, \dotsc,\epsilon(k-1)!\binom{2k}{k}\biggr)
=2\epsilon(k-1)!\binom{2k-1+2\epsilon}{k-1},
\end{equation}
where the generalized binomial coefficient $\binom{z}{k}$ for $z\in\mathbb{C}$ and $k\in\mathbb{N}_0$ is defined by
\begin{equation}\label{falling-binomial-eq}
\binom{z}{k}=
\begin{dcases}
\frac{(-1)^k(-z)_k}{k!}, & k\ge0;\\
0, & k<0.
\end{dcases}
\end{equation}
\end{thm}

\subsection{Alternative proofs of Xu's concise formulas}

In the recent paper~\cite{Filomat-504.tex}, by establishing the formulas
\begin{multline*}
\bell_{2m,k}\biggl(0,\frac{1}{3},0,\frac{1}{5},\dotsc,\frac{1+(-1)^{2m-k+1}}{2}\frac1{2m-k+2}\biggr)\\
=(-1)^k\frac{2^{2m}}{k!}\Biggl[\sum_{\ell=0}^k(-1)^{\ell}\binom{k}{\ell} \frac{T(2m+\ell,\ell)}{\binom{2m+\ell}{\ell}}\Biggr]
\end{multline*}
and
\begin{equation*}
\bell_{2m+1,k}\biggl(0,\frac{1}{3},0,\frac{1}{5},\dotsc,\frac{1+(-1)^{2m-k+2}}{2}\frac1{2m-k+3}\biggr)=0
\end{equation*}
for $m\ge k\in\mathbb{N}_0$, where $\bell_{m,k}$ is defined by~\eqref{partial-Bell-Polyn-Dfn}, verifying the identity
\begin{equation*}
 \sum_{m=0}^j\frac{(-1)^m}{2^m}\binom{j}{m}\sum_{q=0}^m\binom{m}{q}(2q-m)^k
 =\frac{(-1)^j}{2^j} \sum_{\ell=0}^{2j}(-1)^{\ell}\binom{2j}{\ell}(j-\ell)^k
\end{equation*}
for $k\in\mathbb{N}$ and $j\in\mathbb{N}_0$, and making use of the well-known Fa\`a di Bruno formula, the author proved the formulas
\begin{align*}
B_{2k}^{(\epsilon)}\biggl(\frac{\epsilon}{2}\biggr)
&=\sum_{\ell=1}^{2k}\frac{(\epsilon)_\ell}{\ell!} \sum_{j=1}^\ell(-1)^j\binom{\ell}{j} \frac{T(2k+j,j)}{\binom{2k+j}{j}},\\
E_{2k}^{(\epsilon)}\biggl(\frac{\epsilon}{2}\biggr)
&=\frac{1}{4^k} \sum_{\ell=0}^{2k} \frac{(\epsilon)_\ell}{\ell!} \sum_{m=0}^\ell\frac{(-1)^m}{2^m}\binom{\ell}{m}
\sum_{q=0}^m\binom{m}{q}(2q-m)^{2k},
\end{align*}
and
\begin{equation*}
\epsilon(k-1)!\binom{2k-1+\epsilon}{k-1}
=\sum_{\ell=0}^{k}(\epsilon)_{k-\ell} \binom{k+\ell-1}{2\ell} 2^{\ell}(2\ell-1)!!
\end{equation*}
for $k\in\mathbb{N}$ and $\epsilon\in\mathbb{R}$, and then deduced Theorems~\ref{Xu-thm1} through~\ref{Xu-thm3}.

\subsection{The aim of this paper}
The main aim of this paper is to present concise and elegant proofs of Theorems~\ref{Theorem2.1Marian-Gencev-Results-2024-Bell} through~\ref{Theorem2.5Marian-Gencev-Results-2024-Bell} and Theorems~\ref{Xu-thm1} through~\ref{Xu-thm3}. In particular, the proofs of Theorems~\ref{Xu-thm1} through~\ref{Xu-thm3} given below are especially streamlined and elegant.

\section{Concise proofs}

We now proceed to give concise and elegant proofs of Theorems~\ref{Theorem2.1Marian-Gencev-Results-2024-Bell} through~\ref{Theorem2.5Marian-Gencev-Results-2024-Bell}, especially Theorems~\ref{Xu-thm1} through~\ref{Xu-thm3}.

\subsection{Concise proof of Theorem~\ref{Theorem2.1Marian-Gencev-Results-2024-Bell}}\label{proof-theorem5}
Making use of the generating function~\eqref{exp-complete-Bell}, the formula~\eqref{Gencev-ID1-2024-Bell} for $\epsilon=1$ can be reformulated as
\begin{align*}
\exp\Biggl(\frac{1}{2}\sum_{k=1}^{\infty}\frac{B_{2k}}{k}\frac{z^k}{(2k)!}\Biggr)
&=\sum_{k=0}^{\infty}\bell_k\biggl(\frac{1}{2}\frac{B_{2}}{2!}, \frac{1}{2}\frac{B_{4}}{4!}, \dotsc,\frac{1}{2}\frac{(k-1)!B_{2k}}{(2k)!}\biggr)\frac{z^k}{k!}\\
&=\sum_{k=0}^{\infty}\frac{z^k}{(4k)!!(2k+1)}.
\end{align*}
\par
Since
\begin{equation*}
\sinh x=\sum_{k=0}^{\infty}\frac{x^{2k+1}}{(2k+1)!},
\end{equation*}
which can be found in~\cite[p.~42]{Gradshteyn-Ryzhik-Table-8th}, it is easy to see that
\begin{equation*}
\sum_{k=0}^{\infty}\frac{z^k}{(4k)!!(2k+1)}
=\sum_{k=0}^{\infty}\frac{z^k}{2^{2k}(2k)!(2k+1)}
=\sum_{k=0}^{\infty}\frac{(\sqrt{z}\,/2)^{2k}}{(2k+1)!}
=\frac{\sinh(\sqrt{z}\,/2)}{\sqrt{z}\,/2}.
\end{equation*}
Accordingly, it suffices to show
\begin{equation}\label{suffice-cond}
\exp\Biggl(\frac{1}{2}\sum_{k=1}^{\infty}\frac{B_{2k}}{k}\frac{z^k}{(2k)!}\Biggr)
=\frac{\sinh(\sqrt{z}\,/2)}{\sqrt{z}\,/2},
\end{equation}
that is,
\begin{equation}\label{midle-eq}
\frac{1}{2}\sum_{k=1}^{\infty}\frac{B_{2k}}{(2k)!}\frac{z^k}{k}
=\ln\frac{\sinh(\sqrt{z}\,/2)}{\sqrt{z}\,/2}.
\end{equation}
\par
In~\cite[p.~55, Entry~1.518]{Gradshteyn-Ryzhik-Table-8th}, we find the power series expansion
\begin{equation}\label{p.55Item1.518}
\ln\sin x=\ln x+\sum_{k=1}^\infty(-1)^k\frac{2^{2k-1}}{k}B_{2k}\frac{x^{2k}}{(2k)!}, \quad 0<x<\pi.
\end{equation}
The series expansion~\eqref{p.55Item1.518} can be rewritten as
\begin{equation}\label{p.55Item1.518-rewr}
\ln\frac{\sin x}{x}=-\frac{1}{2}\sum_{k=1}^{\infty}\frac{|B_{2k}|}{k}\frac{(2x)^{2k}}{(2k)!}, \quad |x|<\pi.
\end{equation}
From~\eqref{p.55Item1.518-rewr} and the relation $\sinh x=\frac{\sin(\ti x)}{\ti}$ in~\cite[p.~28]{Gradshteyn-Ryzhik-Table-8th}, it follows that
\begin{equation}\label{p.55Item1.518-re}
\ln\frac{\sinh x}{x}=\frac{1}{2}\sum_{k=1}^{\infty}(-1)^{k+1}\frac{|B_{2k}|}{k}\frac{(2x)^{2k}}{(2k)!}, \quad |x|<\pi.
\end{equation}
Replacing $x$ by $\frac{\sqrt{z}\,}{2}$ in~\eqref{p.55Item1.518-re} leads to~\eqref{midle-eq}. Thus, the formula~\eqref{Gencev-ID1-2024-Bell} for $\epsilon=1$ is proved.
\par
Making use of the generating function~\eqref{exp-complete-Bell}, the formula~\eqref{Gencev-ID1-2024-Bell} for $\epsilon=-1$ can be reformulated as
\begin{align*}
\exp\Biggl(-\frac{1}{2}\sum_{k=1}^{\infty}\frac{B_{2k}}{k}\frac{z^k}{(2k)!}\Biggr)
&=\sum_{k=0}^{\infty}\bell_k\biggl(-\frac{1}{2}\frac{B_{2}}{2!}, -\frac{1}{2}\frac{B_{4}}{4!}, \dotsc,-\frac{1}{2}\frac{(k-1)!B_{2k}}{(2k)!}\biggr)\frac{z^k}{k!}\\
&=\sum_{k=0}^{\infty}\frac{2-2^{2k}}{(4k)!!}B_{2k}z^k.
\end{align*}
\par
In~\cite[p.~42]{Gradshteyn-Ryzhik-Table-8th}, we find the series expansion
\begin{equation*}
\csc z=\frac1z+\frac{1}{z}\sum_{k=1}^{\infty}\bigl(2^{2k}-2\bigr)|B_{2k}|\frac{z^{2k}}{(2k)!}, \quad z^2<\pi^2,
\end{equation*}
that is,
\begin{equation*}
\frac{z}{\sin z}=1+\sum_{k=1}^{\infty}\bigl(2^{2k}-2\bigr)|B_{2k}|\frac{z^{2k}}{(2k)!}, \quad z^2<\pi^2.
\end{equation*}
Accordingly, we obtain
\begin{equation*}
\sum_{k=0}^{\infty}\frac{2-2^{2k}}{(4k)!!}B_{2k}z^k
=\sum_{k=0}^{\infty}\bigl(2-2^{2k}\bigr)B_{2k}\frac{(\sqrt{z}\,/2)^{2k}}{(2k)!}
=\frac{\sqrt{-z}\,/2}{\sin (\sqrt{-z}\,/2)}.
\end{equation*}
Consequently, it is sufficient to show
\begin{equation*}
\exp\Biggl(-\frac{1}{2}\sum_{k=1}^{\infty}\frac{B_{2k}}{k}\frac{z^k}{(2k)!}\Biggr)
=\frac{\sqrt{-z}\,/2}{\sin (\sqrt{-z}\,/2)},
\end{equation*}
that is,
\begin{equation*}
-\frac{1}{2}\sum_{k=1}^{\infty}\frac{B_{2k}}{k}\frac{z^k}{(2k)!}
=\ln\frac{\sqrt{-z}\,/2}{\sin (\sqrt{-z}\,/2)}.
\end{equation*}
This follows from replacing $x$ by $\frac{\sqrt{-z}\,}2$ in the series expansion~\eqref{p.55Item1.518-rewr}. Thus, the formula~\eqref{Gencev-ID1-2024-Bell} for $\epsilon=-1$ is proved.
\par
The formula~\eqref{Gencev-ID1-2024-Bell} for $\epsilon=-1$ can be more simply proved as follows.
\par
It is clear that
\begin{align*}
\exp\Biggl(-\frac{1}{2}\sum_{k=1}^{\infty}\frac{B_{2k}}{k}\frac{z^k}{(2k)!}\Biggr)
&=\frac{1}{\exp\bigl(\frac{1}{2}\sum_{k=1}^{\infty}\frac{B_{2k}}{k}\frac{z^k}{(2k)!}\bigr)}\\
&=\frac{\sqrt{z}\,/2}{\sinh(\sqrt{z}\,/2)}\\
&=\sum_{k=0}^{\infty}\frac{2-2^{2k}}{(4k)!!}B_{2k}z^k,
\end{align*}
where we used~\eqref{suffice-cond} and the series expansion
\begin{equation*}
\operatorname{csch} x=\frac{1}{\sinh x}=\frac{1}{x}-\frac{1}{x}\sum_{k=1}^{\infty}\bigl(2^{2k}-1\bigr)B_{2k}\frac{x^{2k}}{(2k)!}, \quad |x|<\pi
\end{equation*}
which can be found in~\cite[p.~42]{Gradshteyn-Ryzhik-Table-8th}. In view of the generating function~\eqref{exp-complete-Bell}, the formula~\eqref{Gencev-ID1-2024-Bell} for $\epsilon=-1$ is thus proved.
The proof of Theorem~\ref{Theorem2.1Marian-Gencev-Results-2024-Bell} is complete.

\subsection{Concise proof of Theorem~\ref{Theorem2.2Marian-Gencev-Results-2024-Bell}}\label{section-thm6p}
Making use of the generating function~\eqref{exp-complete-Bell}, the formula~\eqref{Gencev-ID2-2024-Bell} for $\epsilon=1$ can be reformulated as
\begin{gather*}
\exp\Biggl(\sum_{k=1}^{\infty}\frac{2^{2k}-1}{2}\frac{B_{2k}}{(2k)!}\frac{z^k}{k}\Biggr)
=\sum_{k=0}^{\infty}\bell_k\biggl(\frac{3}{2}\frac{B_{2}}{2!}, \frac{15}{2}\frac{B_{4}}{4!}, \dotsc,\frac{2^{2k}-1}{2}\frac{(k-1)!B_{2k}}{(2k)!}\biggr)\frac{z^k}{k!}\\
=\sum_{k=0}^{\infty}\frac{z^k}{(4k)!!}
=\sum_{k=0}^{\infty}\frac{(\sqrt{z}\,/2)^{2k}}{(2k)!}
=\cos\frac{\sqrt{-z}\,}{2}.
\end{gather*}
Accordingly, it is sufficient to show
\begin{equation}\label{cos-suff}
\exp\Biggl(\sum_{k=1}^{\infty}\frac{2^{2k}-1}{2}\frac{B_{2k}}{(2k)!}\frac{z^k}{k}\Biggr)=\cos\frac{\sqrt{-z}\,}{2},
\end{equation}
that is,
\begin{equation*}
\sum_{k=1}^{\infty}\frac{2^{2k}-1}{2}\frac{B_{2k}}{(2k)!}\frac{z^k}{k}=\ln\cos\frac{\sqrt{-z}\,}{2}.
\end{equation*}
This follows from replacing $x$ by $\frac{\sqrt{-z}\,}{2}$ in the series expansion
\begin{equation*}
\ln\cos x=-\sum_{k=1}^{\infty}\frac{2^{2k-1}(2^{2k}-1)|B_{2k}|}{k}\frac{x^{2k}}{(2k)!}, \quad |x|<\frac{\pi}{2},
\end{equation*}
which can be found in~\cite[p.~55]{Gradshteyn-Ryzhik-Table-8th}. The formula~\eqref{Gencev-ID2-2024-Bell} for $\epsilon=1$ is thus proved.
\par
Making use of the generating function~\eqref{exp-complete-Bell}, the formula~\eqref{Gencev-ID2-2024-Bell} for $\epsilon=-1$ can be reformulated as
\begin{gather*}
\exp\Biggl(\sum_{k=1}^{\infty}-\frac{2^{2k}-1}{2}\frac{B_{2k}}{(2k)!}\frac{z^k}{k}\Biggr)\\
=\sum_{k=0}^{\infty}\bell_k\biggl(-\frac{3}{2}\frac{B_{2}}{2!}, -\frac{15}{2}\frac{B_{4}}{4!}, \dotsc,-\frac{2^{2k}-1}{2}\frac{(k-1)!B_{2k}}{(2k)!}\biggr)\frac{z^k}{k!}\\
=\sum_{k=0}^{\infty}\frac{k!}{(4k)!!}E_{2k}\frac{z^k}{k!}
=\sum_{k=0}^{\infty}\frac{E_{2k}}{(2k)!}\biggl(\frac{\sqrt{z}\,}2\biggr)^{2k}
=\sec\frac{\sqrt{-z}\,}2,
\end{gather*}
where we employed positivity $(-1)^kE_{2k}>0$ and the series expansion
\begin{equation}\label{sec-ser-expan}
\sec x=\sum_{k=0}^{\infty}|E_{2k}|\frac{x^{2k}}{(2k)!}, \quad |x|<\frac{\pi}{2},
\end{equation}
which can be found in~\cite[p.~805]{abram} and~\cite[p.~42]{Gradshteyn-Ryzhik-Table-8th}, respectively.
\par
The formula~\eqref{Gencev-ID2-2024-Bell} for $\epsilon=-1$ can be more simply and concisely proved as follows.
\par
It is obvious that
\begin{gather*}
\exp\Biggl(\sum_{k=1}^{\infty}-\frac{2^{2k}-1}{2}\frac{B_{2k}}{(2k)!}\frac{z^k}{k}\Biggr)
=\frac{1}{\exp\bigl(\sum_{k=1}^{\infty}\frac{2^{2k}-1}{2}\frac{B_{2k}}{(2k)!}\frac{z^k}{k}\bigr)}\\
=\sec\frac{\sqrt{-z}\,}{2}
=\sum_{k=0}^{\infty}\frac{E_{2k}}{(2k)!}\biggl(\frac{\sqrt{z}\,}2\biggr)^{2k}
=\sum_{k=0}^{\infty}\frac{k!}{(4k)!!}E_{2k}\frac{z^k}{k!},
\end{gather*}
where we used the equation~\eqref{cos-suff} and the series expansion~\eqref{sec-ser-expan}. Therefore, the formula~\eqref{Gencev-ID2-2024-Bell} for $\epsilon=-1$ is thus proved.
The proof of Theorem~\ref{Theorem2.2Marian-Gencev-Results-2024-Bell} is thus complete.

\subsection{Concise proof of Theorem~\ref{Theorem2.5Marian-Gencev-Results-2024-Bell}}\label{catalan-sect}
Making use of the generating function~\eqref{exp-complete-Bell}, the formula~\eqref{Gencev-ID3-2024-Bell} for $\epsilon=1$ can be rewritten as
\begin{equation}
\begin{gathered}\label{exp-catanlan}
\exp\Biggl(\frac{1}{2}\sum_{k=1}^{\infty}\binom{2k}{k}\frac{z^k}{k}\Biggr)
=\sum_{k=0}^{\infty}\bell_k\biggl(1, 3, 20, \dotsc, \frac{(k-1)!}{2}\binom{2k}{k}\biggr)\frac{z^k}{k!}\\
=\sum_{k=0}^{\infty}C_kz^k
=\frac2{1+\sqrt{1-4z}\,},
\end{gathered}
\end{equation}
where we used the generating function~\eqref{Catablan-gen-F}. Accordingly, it is sufficient to show
\begin{equation*}
\frac{1}{2}\sum_{k=1}^{\infty}\binom{2k}{k}\frac{z^k}{k}
=\ln\frac2{1+\sqrt{1-4z}\,}.
\end{equation*}
This follows from
\begin{gather*}
\ln\frac2{1+\sqrt{1-4z}\,}=\ln2-\ln\bigl(1+\sqrt{1-4z}\,\bigr)
=\sum_{k=1}^{\infty}(-1)^k\frac{(2k-1)!}{2^{2k}(k!)^2}\bigl(2\sqrt{-z}\,\bigr)^{2k}\\
=\sum_{k=1}^{\infty}(-1)^k\frac{(2k-1)!}{2^{2k}(k!)^2}2^{2k}(-1)^kz^k
=\frac{1}{2}\sum_{k=1}^{\infty}\binom{2k}{k}\frac{z^k}{k},
\end{gather*}
where we used the series expansion
\begin{equation*}
\ln\Bigl(1+\sqrt{1+x^2}\,\Bigr)=\ln 2-\sum_{k=1}^{\infty}(-1)^k\frac{(2k-1)!}{2^{2k}(k!)^2}x^{2k}, \quad x^2\le1,
\end{equation*}
which can be found in~\cite[p.~54]{Gradshteyn-Ryzhik-Table-8th}, and substituting $2\sqrt{-z}\,$ for $x$. The formula~\eqref{Gencev-ID3-2024-Bell} for $\epsilon=1$ is thus proved.
\par
For the case for $\epsilon=-1$, the formula~\eqref{Gencev-ID3-2024-Bell} can be rearranged as
\begin{align*}
\exp\Biggl(-\sum_{k=1}^{\infty}\frac{1}{2}\binom{2k}{k}\frac{z^k}{k}\Biggr)
&=\sum_{k=0}^{\infty}\bell_k\biggl(-1, -3, -20, \dotsc,-\frac{(k-1)!}{2}\binom{2k}{k}\biggr)\frac{z^k}{k!}\\
&=-\sum_{k=1}^{\infty}C_{k-1}z^k,
\end{align*}
that is, in view of~\eqref{exp-catanlan},
\begin{equation*}
\frac{1}{\exp\bigl(\sum_{k=1}^{\infty}\frac{1}{2}\binom{2k}{k}\frac{z^k}{k}\bigr)}
=\frac{1+\sqrt{1-4z}\,}2
=-\sum_{k=1}^{\infty}C_{k-1}z^k.
\end{equation*}
This follows from applying the binomial expansion
\begin{equation*}
(1+x)^q=\sum_{k=0}^{\infty}\binom{q}{k}x^k, \quad q\in\mathbb{R},
\end{equation*}
which can be found in~\cite[p.~25]{Gradshteyn-Ryzhik-Table-8th}, to $\sqrt{1-4z}\,$ by taking $q=\frac{1}{2}$, replacing $x$ by $-4z$, and rearranging carefully, where the generalized binomial coefficient $\binom{q}{k}$ for $q\in\mathbb{C}$ and $k\in\mathbb{N}_0$ is defined by~\eqref{falling-binomial-eq}. The formula~\eqref{Gencev-ID3-2024-Bell} for $\epsilon=-1$ is thus proved. The proof of Theorem~\ref{Theorem2.5Marian-Gencev-Results-2024-Bell} is complete.

\subsection{Concise proof of Theorem~\ref{Xu-thm1}}\label{concise-proofXu-thm1-sec}
From the concise proof of Theorem~\ref{Theorem2.1Marian-Gencev-Results-2024-Bell} in Section~\ref{proof-theorem5}, we conclude that
\begin{equation}\label{exp-bernou-poly}
\exp\Biggl(\epsilon\sum_{k=1}^{\infty}\frac{(k-1)!B_{2k}}{(2k)!}\frac{z^k}{k!}\Biggr)
=\Biggl(\frac{\sin\frac{\sqrt{-z}\,}{2}}{\frac{\sqrt{-z}\,}{2}}\Biggr)^{2\epsilon}, \quad -4\pi^2<z<0.
\end{equation}
It is known~\cite[p.~28]{Gradshteyn-Ryzhik-Table-8th} that
\begin{equation*}
\sin x=\frac{\te^{\ti x}-\te^{-\ti x}}{2\ti}.
\end{equation*}
Then, in light of the generating function~\eqref{Genaralized-BPolyn}, we acquire
\begin{align*}
\biggl(\frac{\sin z}{z}\biggr)^{2\epsilon}&=\biggl(\frac{\te^{\ti z}-\te^{-\ti z}}{2\ti z}\biggr)^{2\epsilon}\\
&=\biggl(\frac{2\ti z}{\te^{2\ti z}-1}\biggr)^{-2\epsilon}\te^{2\ti z(-\epsilon)}\\
&=\sum_{k=0}^{\infty}B_k^{(-2\epsilon)}(-\epsilon)\frac{(2\ti z)^k}{k!}\\
&=\sum_{k=0}^{\infty}B_k^{(-2\epsilon)}(-\epsilon)\frac{(2z)^k\ti^k}{k!}\\
&=\sum_{k=0}^{\infty}B_{2k}^{(-2\epsilon)}(-\epsilon)\frac{(2z)^{2k}\ti^{2k}}{(2k)!} +\sum_{k=0}^{\infty}B_{2k+1}^{(-2\epsilon)}(-\epsilon)\frac{(2z)^{2k+1}\ti^{2k+1}}{(2k+1)!}\\
&=\sum_{k=0}^{\infty}(-1)^kB_{2k}^{(-2\epsilon)}(-\epsilon)\frac{(2z)^{2k}}{(2k)!} +\ti\sum_{k=0}^{\infty}(-1)^kB_{2k+1}^{(-2\epsilon)}(-\epsilon)\frac{(2z)^{2k+1}}{(2k+1)!}\\
&=\sum_{k=0}^{\infty}(-1)^kB_{2k}^{(-2\epsilon)}(-\epsilon)\frac{(2z)^{2k}}{(2k)!}
\end{align*}
for $|z|<\pi$. Therefore, it follows that
\begin{equation*}
\Biggl(\frac{\sin\frac{\sqrt{-z}\,}{2}}{\frac{\sqrt{-z}\,}{2}}\Biggr)^{2\epsilon}
=\sum_{k=0}^{\infty}(-1)^kB_{2k}^{(-2\epsilon)}(-\epsilon)\frac{\bigl(\sqrt{-z}\,\bigr)^{2k}}{(2k)!}
=\sum_{k=0}^{\infty}B_{2k}^{(-2\epsilon)}(-\epsilon)\frac{z^k}{(2k)!}.
\end{equation*}
Combining this series expansion with~\eqref{exp-bernou-poly} yields
\begin{equation*}
\exp\Biggl(\epsilon\sum_{k=1}^{\infty}\frac{(k-1)!B_{2k}}{(2k)!}\frac{z^k}{k!}\Biggr)
=\sum_{k=0}^{\infty}B_{2k}^{(-2\epsilon)}(-\epsilon)\frac{k!}{(2k)!}\frac{z^k}{k!}.
\end{equation*}
Comparing this with the generating function~\eqref{exp-complete-Bell}, we arrive at the formula~\eqref{Xu-Eq1}.
The proof of Theorem~\ref{Xu-thm1} is complete.

\subsection{Concise proof of Theorem~\ref{Xu-thm2}}\label{concise-proofXu-thm2-sec}
Basing on the proof of Theorem~\ref{Theorem2.2Marian-Gencev-Results-2024-Bell} in Section~\ref{section-thm6p}, it follows that
\begin{equation}\label{cos-epsilon-exp}
\exp\Biggl(\epsilon\sum_{k=1}^{\infty}\bigl(2^{2k}-1\bigr)\frac{(k-1)!B_{2k}}{(2k)!}\frac{z^k}{k!}\Biggr)
=\biggl(\cos\frac{\sqrt{-z}\,}{2}\biggr)^{2\epsilon}.
\end{equation}
\par
Making use of the exponential representation
\begin{equation*}
\cos x=\frac{\te^{\ti x}+\te^{-\ti x}}{2}
\end{equation*}
in~\cite[p.~28]{Gradshteyn-Ryzhik-Table-8th} and the generating function~\eqref{Gen-Euler-Polyn}, we derive
\begin{align*}
(\cos z)^{2\epsilon}&=\biggl(\frac{\te^{\ti z}+\te^{-\ti z}}{2}\biggr)^{2\epsilon}\\
&=\biggl(\frac{2}{\te^{2\ti z}+1}\biggr)^{-2\epsilon}\te^{-2\epsilon\ti z}\\
&=\sum_{k=0}^{\infty}E_k^{(-2\epsilon)}(-\epsilon)\frac{(2\ti z)^k}{k!}\\
&=\sum_{k=0}^{\infty}E_{2k}^{(-2\epsilon)}(-\epsilon)\frac{(2\ti z)^{2k}}{(2k)!} +\sum_{k=0}^{\infty}E_{2k+1}^{(-2\epsilon)}(-\epsilon)\frac{(2\ti z)^{2k+1}}{(2k+1)!}\\
&=\sum_{k=0}^{\infty}(-1)^kE_{2k}^{(-2\epsilon)}(-\epsilon)\frac{(2z)^{2k}}{(2k)!} +\ti\sum_{k=0}^{\infty}(-1)^kE_{2k+1}^{(-2\epsilon)}(-\epsilon)\frac{(2z)^{2k+1}}{(2k+1)!}\\
&=\sum_{k=0}^{\infty}(-1)^kE_{2k}^{(-2\epsilon)}(-\epsilon)\frac{(2z)^{2k}}{(2k)!}.
\end{align*}
Accordingly, we have
\begin{equation*}
\biggl(\cos\frac{\sqrt{-z}\,}{2}\biggr)^{2\epsilon}
=\sum_{k=0}^{\infty}E_{2k}^{(-2\epsilon)}(-\epsilon)\frac{z^{k}}{(2k)!}.
\end{equation*}
Substituting this into~\eqref{cos-epsilon-exp} gives
\begin{equation*}
\exp\Biggl(\epsilon\sum_{k=1}^{\infty}\bigl(2^{2k}-1\bigr)\frac{(k-1)!B_{2k}}{(2k)!}\frac{z^k}{k!}\Biggr)
=\sum_{k=0}^{\infty}\frac{k!}{(2k)!}E_{2k}^{(-2\epsilon)}(-\epsilon)\frac{z^{k}}{k!}.
\end{equation*}
Comparing this with the generating function~\eqref{exp-complete-Bell}, we arrive at the formula~\eqref{Xu-Eq2}.
The proof of Theorem~\ref{Xu-thm2} is complete.

\subsection{Concise proof of Theorem~\ref{Xu-thm3}}
Basing the concise proof of Theorem~\ref{Theorem2.5Marian-Gencev-Results-2024-Bell} in Section~\ref{catalan-sect}, we obtain
\begin{equation}\label{exp=2epsilon=pow}
\exp\Biggl(\epsilon\sum_{k=1}^{\infty}\binom{2k}{k}\frac{z^k}{k}\Biggr)
=\biggl(\frac2{1+\sqrt{1-4z}\,}\biggr)^{2\epsilon}.
\end{equation}
\par
In~\cite[p.~25, Entry~1.114]{Gradshteyn-Ryzhik-Table-8th}, we find that
\begin{equation*}
\bigl(1+\sqrt{1+x}\,\bigr)^q
=2^q\biggl[1+\frac{q}{1!}\biggl(\frac{x}{4}\biggr) +\frac{q(q-3)}{2!}\biggl(\frac{x}{4}\biggr)^2 +\frac{q(q-4)(q-5)}{3!}\biggl(\frac{x}{4}\biggr)^3+\dotsm\biggr]
\end{equation*}
for $x^2<1$ and $q\in\mathbb{R}$. This series expansion can be reformulated as
\begin{equation}\label{Item1.114p.25GRT-8th-Gen}
\biggl(\frac{1+\sqrt{1+x}\,}{2}\biggr)^q
=1+q\sum_{k=1}^{\infty} \binom{q-k-1}{k-1}\frac{1}{k}\biggl(\frac{x}{4}\biggr)^k
\end{equation}
for $|x|<1$ and $q\in\mathbb{R}$, which will be proved in Section~\ref{Appendix} below. As a result, we acquire
\begin{align*}
\biggl(\frac2{1+\sqrt{1-4z}\,}\biggr)^{2\epsilon}
&=\biggl(\frac{1+\sqrt{1-4z}\,}2\biggr)^{-2\epsilon}\\
&=1-2\epsilon\sum_{k=1}^{\infty} (-1)^k(k-1)!\binom{-2\epsilon-k-1}{k-1}\frac{z^k}{k!}\\
&=1+2\epsilon\sum_{k=1}^{\infty} (k-1)!\binom{2\epsilon+2k-1}{k-1}\frac{z^k}{k!},
\end{align*}
where we used the transform
\begin{equation*}
\binom{x}{k}=(-1)^k\binom{-x+k-1}{k}, \quad (x,k)\in\mathbb{C}\times\mathbb{N}_0,
\end{equation*}
which can be found in~\cite[p.~3, (1.6)]{Quaintance-Gould-Stirling-B}. Substituting this result into~\eqref{exp=2epsilon=pow} leads to
\begin{equation*}
\exp\Biggl(\epsilon\sum_{k=1}^{\infty}\binom{2k}{k}\frac{z^k}{k}\Biggr)
=1+2\epsilon\sum_{k=1}^{\infty} (k-1)!\binom{2\epsilon+2k-1}{k-1}\frac{z^k}{k!}.
\end{equation*}
Comparing this with the generating function~\eqref{exp-complete-Bell}, we arrive at the formula~\eqref{Xu-Eq3}.
The proof of Theorem~\ref{Xu-thm3} is complete.

\section{Remarks}
In this section, we give several remarks about some stuffs related to our concise and elegant proofs.

\begin{rem}
From the concise proof of Theorem~\ref{Xu-thm1} in Section~\ref{concise-proofXu-thm1-sec}, we conclude that the generalized Bernoulli polynomials satisfy
$$
B_{2k+1}^{(2\epsilon)}(\epsilon)=0
$$
for $\epsilon\in\mathbb{C}$ and $k\in\mathbb{N}_0$ and that the series expansion
\begin{equation}\label{sinc-ser-expan}
\biggl(\frac{\sin z}{z}\biggr)^{r}=\sum_{k=0}^{\infty}(-1)^kB_{2k}^{(-r)}\biggl(-\frac{r}{2}\biggr)\frac{(2z)^{2k}}{(2k)!}
\end{equation}
is true for $r\in\mathbb{R}$. The form and its proof of the series expansion~\eqref{sinc-ser-expan} are more concise than those of the series expansion
\begin{equation*}
\biggl(\frac{\sin z}{z}\biggr)^r=1+\sum_{q=1}^{\infty}(-1)^q\Biggl[\sum_{\ell=1}^{2q}\frac{(-r)_\ell}{\ell!}
\sum_{k=1}^\ell(-1)^k\binom{\ell}{k} \frac{T(2q+k,k)}{\binom{2q+k}{k}}\Biggr]\frac{(2z)^{2q}}{(2q)!}
\end{equation*}
for $r\in\mathbb{R}$, which can be found in~\cite[Example~3.2]{integer2Bell2real.tex} and~\cite[Theorem~11]{AADM-3223-2023.tex}, where the rising factorial $(z)_k$ for $z\in\mathbb{C}$ and $k\in\mathbb{N}_0$ is defined by~\eqref{rising-Factorial}.
\par
The series expansion~\eqref{sinc-ser-expan} recovers~\cite[p.~340, Eq.~(50.8.10)]{Hansen-B-1975} which reads that
\begin{equation*}
\sum_{k=0}^{\infty}(-1)^k\frac{a^{2k}}{(2k)!}B_{2k}^{(-m)}\biggl(-\frac{m}{2}\biggr)=\biggl(\frac{2}{a}\sin\frac{a}{2}\biggr)^m.
\end{equation*}
\end{rem}

\begin{rem}
From the concise proof of Theorem~\ref{Xu-thm2} in Section~\ref{concise-proofXu-thm2-sec}, we conclude that the generalized Euler polynomials satisfy
$$
E_{2k+1}^{(2\epsilon)}(\epsilon)=0
$$
for $\epsilon\in\mathbb{C}$ and $k\in\mathbb{N}_0$ and that the series expansion
\begin{equation}\label{cosine-power-series}
(\cos z)^r=\sum_{k=0}^{\infty}(-1)^kE_{2k}^{(-r)}\biggl(-\frac{r}{2}\biggr)\frac{(2z)^{2k}}{(2k)!}
\end{equation}
is true for $r\in\mathbb{R}$. The form and its proof of the series expansion~\eqref{cosine-power-series} are more concise than those of the series expansion
\begin{equation*}
(\cos z)^{r}
=\sum_{k=0}^{\infty}(-1)^k \Biggl[\sum_{\ell=0}^{2k} \frac{(-r)_\ell}{\ell!} \sum_{m=0}^\ell\frac{(-1)^m}{2^m}\binom{\ell}{m}
\sum_{q=0}^m\binom{m}{q}(2q-m)^{2k}\Biggr]\frac{z^{2k}}{(2k)!}
\end{equation*}
for $r\in\mathbb{R}$, which was established in~\cite[pp.~10--11, Eq.~(32)]{RIMA-D-23-01264.tex}.
\end{rem}

\begin{rem}
Comparing Theorems~\ref{Theorem2.1Marian-Gencev-Results-2024-Bell} and~\ref{Xu-thm1}, we deduce
\begin{equation*}
B_{2k}^{(-1)}\biggl(-\frac{1}{2}\biggr)=\frac{1}{2^{2k}(2k+1)}
\quad\text{and}\quad
B_{2k}^{(1)}\biggl(\frac{1}{2}\biggr)=\frac{2-2^{2k}}{2^{2k}}B_{2k}
\end{equation*}
for $k\in\mathbb{N}_0$.
\end{rem}

\begin{rem}
Comparing Theorems~\ref{Theorem2.2Marian-Gencev-Results-2024-Bell} and~\ref{Xu-thm2}, we deduce
\begin{equation*}
E_{2k}^{(-1)}\biggl(-\frac{1}{2}\biggr)=\frac{1}{2^{2k}}
\quad\text{and}\quad
E_{2k}^{(1)}\biggl(\frac{1}{2}\biggr)=\frac{1}{2^{2k}}E_{2k}
\end{equation*}
for $k\in\mathbb{N}_0$.
\end{rem}

\section{Appendix}\label{Appendix}
In this section, we give a proof of the series expansion~\eqref{Item1.114p.25GRT-8th-Gen} in details.
\par
In~\cite[Theorem~11.4]{Charalambides-book-2002}, the monograph~\cite{Ion-Chitescu-B-2017}, and~\cite[p.~139, Theorem~C]{Comtet-Combinatorics-74}, we find the famous Fa\`a di Bruno formula
\begin{equation}\label{Bruno-Bell-Polynomial}
\frac{\operatorname{d}^k}{\td z^k}f\circ h(z)=\sum_{j=0}^k f^{(j)}(h(z)) \bell_{k,j}\bigl(h'(z),h''(z),\dotsc,h^{(k-j+1)}(z)\bigr)
\end{equation}
for $k\in\mathbb{N}_0$, where $\bell_{k,j}(z_1,z_2,\dotsc,z_{k-j+1})$ is defined by~\eqref{partial-Bell-Polyn-Dfn}, the function $f$ is $k$-time differentiable, and $h$ is $(k+1)$-time differentiable. Applying~\eqref{Bruno-Bell-Polynomial}, we obtain
\begin{multline*}
\frac{\operatorname{d}^k}{\td x^k}\biggl[\biggl(\frac{1+\sqrt{1+x}\,}{2}\biggr)^q\biggr]
=\sum_{j=1}^{k} (u^q)^{(j)} \bell_{k,j}\bigl(u'(x),u''(x),\dotsc,u^{(k-j+1)}(x)\bigr)\\
\begin{aligned}
&=\sum_{j=1}^{k} \langle q\rangle_j u^{q-j} \bell_{k,j}\biggl(\frac{1}{2} \biggl\langle\frac{1}{2}\biggr\rangle_1(1+x)^{1/2-1}, \frac{1}{2} \biggl\langle\frac{1}{2}\biggr\rangle_2(1+x)^{1/2-2}, \\
&\quad\dotsc, \frac{1}{2} \biggl\langle\frac{1}{2}\biggr\rangle_{k-j+1}(1+x)^{1/2-(k-j+1)}\biggr),
\end{aligned}
\end{multline*}
where $u=u(x)=\frac{1+\sqrt{1+x}\,}{2}$ and the falling factorial $\langle z\rangle_n$ for $z\in\mathbb{C}$ and $n\in\mathbb{N}_0$ is defined by
\begin{equation*}
\langle z\rangle_n=
\prod_{k=0}^{n-1}(z-k)=
\begin{cases}
z(z-1)\dotsm(z-n+1), & n\in\mathbb{N};\\
1,& n=0.
\end{cases}
\end{equation*}
\par
Further making use of the identity
\begin{equation*}
\bell_{k,j}\bigl(abz_1,ab^2z_2,\dotsc,ab^{k-j+1}z_{k-j+1}\bigr) =a^j b^k \bell_{k,j}(z_1,z_2,\dotsc,z_{k-j+1})
\end{equation*}
and the formula
\begin{equation*}
\bell_{k,j}\biggl(\biggl\langle\frac12\biggr\rangle_1, \biggl\langle\frac12\biggr\rangle_2, \dotsc, \biggl\langle\frac12\biggr\rangle_{k-j+1}\biggr)
=(-1)^{k+j}\frac{(2k-2j-1)!!}{2^k}\binom{2k-j-1}{2k-2j}
\end{equation*}
for $k\ge j\in\mathbb{N}_0$, which can be found in~\cite[p.~412]{Charalambides-book-2002}, \cite[p.~135]{Comtet-Combinatorics-74}, \cite[Theorem~4.1]{Special-Bell2Euler.tex}, and the proof of~\cite[Theorem~3.2]{CDM-68111.tex}, respectively, we derive
\begin{align*}
&\quad \frac{\operatorname{d}^k}{\td x^k}\biggl[\biggl(\frac{1+\sqrt{1+x}\,}{2}\biggr)^q\biggr]\\
&=\sum_{j=1}^{k} \langle q\rangle_j u^{q-j}(x)\frac{(1+x)^{j/2-k}}{2^j} \bell_{k,j}\biggl(\biggl\langle\frac{1}{2}\biggr\rangle_1, \biggl\langle\frac{1}{2}\biggr\rangle_2, \dotsc, \biggl\langle\frac{1}{2}\biggr\rangle_{k-j+1}\biggr)\\
&\to\sum_{j=1}^{k} \frac{\langle q\rangle_j}{2^j} \bell_{k,j}\biggl(\biggl\langle\frac{1}{2}\biggr\rangle_1, \biggl\langle\frac{1}{2}\biggr\rangle_2, \dotsc, \biggl\langle\frac{1}{2}\biggr\rangle_{k-j+1}\biggr), \quad x\to0\\
&=\sum_{j=1}^{k} (-1)^{k+j}\frac{\langle q\rangle_j}{2^{k+j}} (2k-2j-1)!!\binom{2k-j-1}{2k-2j}
\end{align*}
for $(q,k)\in\mathbb{R}\times\mathbb{N}$. Accordingly, we arrive at
\begin{equation*}
\biggl(\frac{1+\sqrt{1+x}\,}{2}\biggr)^q
=1+\sum_{k=1}^{\infty} \Biggl[\sum_{j=1}^{k} (-1)^{k+j}\frac{\langle q\rangle_j}{2^{k+j}} (2k-2j-1)!!\binom{2k-j-1}{2k-2j}\Biggr] \frac{x^k}{k!}
\end{equation*}
for $(q,x)\in\mathbb{R}\times(-1,1)$. Hence, in order to prove the series expansion~\eqref{Item1.114p.25GRT-8th-Gen}, it is enough to show
\begin{equation*}
\sum_{j=1}^{k} (-1)^{k+j}\frac{\langle q\rangle_j}{2^{k+j}} (2k-2j-1)!!\binom{2k-j-1}{2k-2j}
=q\binom{q-k-1}{k-1}\frac{(k-1)!}{4^k}
\end{equation*}
for $(q,k)\in\mathbb{R}\times\mathbb{N}$, which is equivalent to
\begin{equation}\label{Zhang-Lim-Qi-equiv}
\sum_{j=0}^{k}(-1)^j\binom{q}{j}\binom{2k-j}{k} =(-1)^{k}\binom{q-k-1}{k}
\end{equation}
for $(q,k)\in\mathbb{R}\times\mathbb{N}_0$. The identity~\eqref{Zhang-Lim-Qi-equiv} follows from the combination of the identity
\begin{equation*}
\sum_{k=0}^{n}(-1)^k\binom{x}{k}\binom{2n-k}{n}=\binom{2n-x}{n}, \quad (x,n)\in\mathbb{C}\times\mathbb{N}_0,
\end{equation*}
see~\cite[p.~65, (3.50)]{Sprugnoli-Gould-2006}, with the equality
\begin{equation*}
\binom{x}{k}=(-1)^k\binom{-x+k-1}{k}, \quad (x,k)\in\mathbb{C}\times\mathbb{N}_0,
\end{equation*}
see~\cite[p.~3, (1.6)]{Quaintance-Gould-Stirling-B}. The proof of the series expansion~\eqref{Item1.114p.25GRT-8th-Gen} is complete.

\section{Declarations}

\paragraph{\bf Authors' Contributions}
All authors contributed equally to the manuscript and read and approved the final manuscript.

\paragraph{\bf Funding}
The author was partially supported by the Natural Science Foundation of Inner Mongolia Autonomous Region (Grant No.~2025QN01041) and by the Youth Project of Hulunbuir City for Basic Research and Applied Basic Research (Grant No.~GH2024020)

\paragraph{\bf Acknowledgements}
Not applicable.

\paragraph{\bf Institutional Review Board Statement}
Not applicable.

\paragraph{\bf Informed Consent Statement}
Not applicable.

\paragraph{\bf Ethical Approval}
The conducted research is not related to either human or animal use.

\paragraph{\bf Availability of Data and Material}
Data sharing is not applicable to this article as no new data were created or analyzed in this study.

\paragraph{\bf Competing Interests}
The authors declare that they have no conflict of competing interests.

\paragraph{\bf Use of AI Tools Declaration}
The authors declare they have not used Artificial Intelligence (AI) tools in the creation of this article.


\begin{thebibliography}{99}

\bibitem{abram}
M. Abramowitz and I. A. Stegun (Eds), \textit{Handbook of Mathematical Functions with Formulas, Graphs, and Mathematical Tables}, National Bureau of Standards, Applied Mathematics Series \textbf{55}, Reprint of the 1972 edition, Dover Publications, Inc., New York, 1992.

\bibitem{MIA-9509.tex}
L.-Y. Bao, C.-Y. He, and F. Qi, \textit{Monotonic sequences and inequalities involving the ratio between two adjacent nonzero Bernoulli numbers}, Math. Inequal. Appl. \textbf{29} (2026), no.~1, 1\nobreakdash--14. DOI: \url{https://doi.org/10.7153/mia-2026-29-01}.

\bibitem{Charalambides-book-2002}
C. A. Charalambides, \textit{Enumerative Combinatorics}, CRC Press Series on Discrete Mathematics and its Applications. Chapman \& Hall/CRC, Boca Raton, FL, 2002.

\bibitem{Ion-Chitescu-B-2017}
I. Chi\c{t}escu, \textit{Around the formula of Fa\`a di Bruno}, \'Editions universitaires europ\'eennes, Mauritius, 2017.

\bibitem{Comtet-Combinatorics-74}
L. Comtet, \emph{Advanced Combinatorics: The Art of Finite and Infinite Expansions}, Revised and Enlarged Edition, D. Reidel Publishing Co., Dordrecht and Boston, 1974. DOI: \url{https://doi.org/10.1007/978-94-010-2196-8}.

\bibitem{Marian-Gencev-Results-2024}
M. Gen\v{c}ev, \emph{Extension of Hoffman's combinatorial identity via specific zeta-like series}, Results Math. \textbf{79} (2024), no.~1, Paper No.~2, 25~pages. DOI: \url{https://doi.org/10.1007/s00025-023-02035-w}.

\bibitem{Gil-Fresan-2020}
J. I. B. Gil and J. Fres\'an, \emph{Multiple Zeta Values: from Numbers to Motives}, Clay Mathematics Proceedings, 2020. URL: \url{https://javier.fresan.perso.math.cnrs.fr/mzv.pdf}.

\bibitem{Gradshteyn-Ryzhik-Table-8th}
I. S. Gradshteyn and I. M. Ryzhik, \emph{Table of Integrals, Series, and Products}, Translated from the Russian, Translation edited and with a preface by Daniel Zwillinger and Victor Moll, Eighth edition, Revised from the seventh edition, Elsevier/Academic Press, Amsterdam, 2015. DOI: \url{https://doi.org/10.1016/B978-0-12-384933-5.00013-8}.

\bibitem{Hansen-B-1975}
E. R. Hansen, \emph{A Table of Series and Products}, Prentice-Hall, Englewood Cliffs, NJ, USA, 1975.

\bibitem{RIMA-D-23-01264.tex}
C.-Y. He and F. Qi, \textit{Reformulations and generalizations of Hoffman's and Gen\v{c}ev's combinatorial identities}, Results Math. \textbf{79} (2024), no.~4, Paper No.~131, 17~pages. DOI: \url{https://doi.org/10.1007/s00025-024-02160-0}.

\bibitem{Hoffman-Pacific-1992}
M. E. Hoffman, \emph{Multiple harmonic series}, Pacific J. Math. \textbf{152} (1992), no.~2, 275\nobreakdash--290. URL: \url{http://projecteuclid.org/euclid.pjm/1102636166}.

\bibitem{Koshy-B-2009}
T. Koshy, \emph{Catalan Numbers with Applications}, Oxford University Press, Oxford, 2009.

\bibitem{Sprugnoli-Gould-2006}
R. Sprugnoli, \textit{Riordan Array Proofs of Identities in Gould's Book}, University of Florence, Italy, 2006.

\bibitem{Filomat-504.tex}
F. Qi, \textit{Alternative proofs of Xu's forms for He--Qi's combinatorial identities reformulating and generalizing Hoffman's and Gen\v{c}ev's combinatorial identities for complete Bell polynomials}, Filomat \textbf{40} (2026), accepted on 10 February 2026. DOI: \url{https://www.researchgate.net/publication/400661711}.

\bibitem{AADM-3251.tex}
F. Qi, \emph{Some properties of the Catalan numbers}, Appl. Anal. Discrete Math. \textbf{19} (2025), no.~1, 176\nobreakdash--184. DOI: \url{https://doi.org/10.2298/AADM240130002Q}.

\bibitem{Special-Bell2Euler.tex}
F. Qi and B.-N. Guo, \textit{Explicit formulas for special values of the Bell polynomials of the second kind and for the Euler numbers and polynomials}, Mediterr. J. Math. \textbf{14} (2017), no.~3, Article~140, 14~pages. DOI: \url{https://doi.org/10.1007/s00009-017-0939-1}.

\bibitem{integer2Bell2real.tex}
F. Qi, G. V. Milovanovi\'c, and D. Lim, \emph{Specific values of partial Bell polynomials and series expansions for real powers of functions and for composite functions}, Filomat \textbf{37} (2023), no.~28, 9469\nobreakdash--9485. DOI: \url{https://doi.org/10.2298/FIL2328469Q}.

\bibitem{CDM-68111.tex}
F. Qi, D.-W. Niu, D. Lim, and B.-N. Guo, \textit{Closed formulas and identities for the Bell polynomials and falling factorials}, Contrib. Discrete Math. \textbf{15} (2020), no.~1, 163\nobreakdash--174. DOI: \url{https://doi.org/10.11575/cdm.v15i1.68111}.

\bibitem{AADM-3223-2023.tex}
F. Qi and P. Taylor, \emph{Series expansions for powers of sinc function and closed-form expressions for specific partial Bell polynomials}, Appl. Anal. Discrete Math. \textbf{18} (2024), no.~1, 92\nobreakdash--115. DOI: \url{https://doi.org/10.2298/AADM230902020Q}.

\bibitem{Quaintance-Gould-Stirling-B}
J. Quaintance and H. W. Gould, \emph{Combinatorial Identities for Stirling Numbers}, the unpublished notes of H. W. Gould, with a foreword by George E. Andrews, World Scientific Publishing Co. Pte. Ltd., Singapore, 2016.

\bibitem{Riordan-B-1968}
J. Riordan, \emph{Combinatorial Identities}, Reprint of the 1968 original, Robert E. Krieger Publishing Co., Huntington, N.Y., 1979.

\bibitem{Stanley-Catalan-2015}
R. P. Stanley, \emph{Catalan Numbers}, Cambridge University Press, New York, 2015. DOI: \url{https://doi.org/10.1017/CBO9781139871495}.

\bibitem{Temme-96-book}
N. M. Temme, \emph{Special Functions: An Introduction to Classical Functions of Mathematical Physics}, A Wiley-Interscience Publication, John Wiley \& Sons, Inc., New York, 1996. DOI: \url{http://dx.doi.org/10.1002/9781118032572}.

\bibitem{Xu-AIMS-Math-2025}
A. Xu, \textit{Concise forms of He--Qi combinatorial identities in terms of generalized Bernoulli polynomials, Euler polynomials, and Catalan numbers}, ResearchGate Preprint (2026). DOI: \url{https://doi.org/10.13140/RG.2.2.14137.63849}.

\bibitem{Don.1992-94}
D. Zagier, \emph{Values of zeta functions and their applications}, First European Congress of Mathematics, Vol. II (Paris, 1992), 497\nobreakdash--512, Progr. Math., 120, Birkh\"auser, Basel, 1994.

\bibitem{Zhao-zeta-multiple}
J. Zhao, \emph{Multiple Zeta Functions, Multiple Polylogarithms and Their Special Values}, Series on Number Theory and its Applications, 12, World Scientific Publishing Co. Pte. Ltd., Hackensack, NJ, 2016. DOI: \url{https://doi.org/10.1142/9634}.

\end{thebibliography}
\end{document}